\documentclass[11pt]{amsart}
\usepackage{graphicx}
\usepackage[centertags]{amsmath}
\usepackage{amsfonts}
\usepackage{amssymb}
\usepackage{amsthm}
\usepackage{newlfont}
\newtheorem{theorem}{Theorem}[section]
\newtheorem{corollary}[theorem]{Corollary}

\newtheorem{prop}[theorem]{Proposition}
\theoremstyle{definition}
\newtheorem{defn}[theorem]{Definition}
\theoremstyle{remark}

\newcommand{\id}{\textrm{id}}
\newcommand{\grad}{\textrm{grad}}

\title[Three-dimensional Riemannian manifolds with circulant structures]{Three-dimensional Riemannian manifolds with circulant structures}

\author{Iva Dokuzova}
\address{Department of Algebra and Geometry, Faculty of Mathematics and Informatics, University of
Plovdiv Paisii Hilendarski, 236 Bulgaria Blvd, 4027 Plovdiv, Bulgaria}
\email{dokuzova@uni-plovdiv.bg}

\author{ D. Razpopov}
\address{Department of Mathematics and Physics, Agricultural University of Plovdiv, 12 Mendeleev Blvd., 4000 Plovdiv, Bulgaria}
\email{razpopov@au-plovdiv.bg}

\author{G. Dzhelepov}
\address{Department of Mathematics and Physics, Agricultural University of Plovdiv, 12 Mendeleev Blvd., 4000 Plovdiv, Bulgaria}
\email{dzhelepov@au-plovdiv.bg}

\begin{document}

\begin{abstract}
We consider a $3$-dimensional Riemannian manifold $M$ with two circulant structures -- a metric $g$ and an endomorphism $q$ whose third power is identity. The structure $q$ is compatible with $g$ such that an isometry is induced in any tangent space of $M$.  We obtain some curvature properties of this manifold $(M, g, q)$ and give an explicit example of such a manifold.
\end{abstract}

\maketitle
\section{Introduction}
\label{intro}

Circulant matrices occur in many areas of the applied mathematics. For instance, they are particulary useful in the Vibration analysis, Linear codes, Geometry, Graph theory, etc. (see \cite{72}, \cite{22}, \cite{92}, \cite{13}). We are motivated to equip differentiable manifolds with additional structures which are represented by circulant matrices.

In differential geometry essential results are associated with the sectional curvatures of some characteristic $2$-planes of the tangent space of the manifolds with additional structures, for example  \cite{132}, \cite{55} and \cite{52}. Another important problem is the obtaining of explicit examples of the constructed manifolds.

The main aim of the present study is on differential geometry of $3$-dimensional Riemannian manifolds equipped with an endomorphism $q$ whose third power is identity. Moreover the metric $g$ and the structure $q$ are represented by circulant matrices.

The paper is organized as follows. In Sect.~\ref{sec:2}, we consider a $3$-dimensional Riemannian manifold $M$ with a circulant metric and a circulant structure $q$ with $q^{3}=\id$, i.e. a manifold $(M, g, q)$. Also we recall necessary facts about such manifolds and about a $q$-\textit{basis} of the tangent space $T_{p}M$, $p\in M$. In Sect.~\ref{sec:3}, we calculate the components of the corresponding curvature tensor $R$ with respect to the Levi-Civita connection of  $g$. In Sect.~\ref{sec:4}, we consider two special properties of $R$ with respect to $q$ and the consequences for some sectional curvatures.  In Sect.~\ref{sec:5}, we obtain an explicit example.

\section{Preliminaries}
\label{sec:2}
Let $M$ be a $3$-dimensional manifold with a Riemannian metric $g$. Let the components of the metric $g$ at an arbitrary point $p(X^{1}, X^{2}, X^{3})\in M$ form the following circulant matrix
\begin{equation}\label{f2}
    (g_{ij})=\begin{pmatrix}
      A & B & B \\
      B & A & B \\
      B & B & A \\
    \end{pmatrix},
\end{equation}
where $A$ and $B$ are smooth functions of $X^{1}$, $X^{2}$, $X^{3}$.

We assume that
\begin{equation}\label{f3}
  A > B > 0.
\end{equation}
Then the conditions to be a positive definite metric $g$ are satisfied:
\begin{equation*}
    A>0,\quad \begin{vmatrix}
      A & B\\
      B & A \\
    \end{vmatrix}=(A-B)(A+B)>0,
\end{equation*}
\begin{equation*}
\begin{vmatrix}
      A & B & B \\
      B & A & B \\
      B & B & A \\
    \end{vmatrix}=(A-B)^{2}(A+2B)>0.
\end{equation*}

Let $q$ be an endomorphism in the tangent space $T_{p}M$, whose coordinate matrix with respect to a basis $\{e_{i}\}$ of $T_{p}M$ is
\begin{equation}\label{f4}
    (q_{i}^{.j})=\begin{pmatrix}
      0 & 1 & 0 \\
      0 & 0 & 1 \\
      1 & 0 & 0 \\
    \end{pmatrix}.
\end{equation}
Then
\begin{equation*}
    q^{3}= \id.
\end{equation*}

We denote by $(M, g, q)$ the manifold $M$ equipped with the metric $g$ and the structure $q$, which are defined by (\ref{f2}) -- (\ref{f4}).

Further, $x, y, z, u$ will stand for arbitrary elements of the algebra on the smooth vector fields on $M$ or vectors in the tangent space $T_{p}M$. The Einstein summation convention is used, the range of the summation indices being always $\{1, 2, 3\}$.

In \cite{1} it is proved that the structure $q$ of the manifold $(M, g, q)$ is an isometry with respect to the metric $g$, i.e.
\begin{equation}\label{2.12}
     g(qx, qy)=g(x,y).
\end{equation}

\begin{defn}
A basis of type $\{x, qx, q^{2}x\}$ of $T_{p}M$ is called a $q$-\textit{basis}. In this case we say that \textit{the vector $x$ induces a $q$-basis of} $T_{p}M$. Similarly, a basis $\{x, qx\}$ of a $2$-plane $\alpha=\{x, qx\}$ is called a $q$-\textit{basis}.
\end{defn}

In \cite{1} it is verified that
\begin{itemize}
\item [(i)]  A vector $x=(x^{1}, x^{2}, x^{3})$ induces a $q$-basis in $T_{p}M$ if and only if
   \begin{equation*}
   3x^{1}x^{2}x^{3}\neq (x^{1})^{3}+(x^{2})^{3}+(x^{3})^{3};
   \end{equation*}
  \item [(ii)] If a vector $x$ induces a $q$-basis of $T_{p}M$ and  $\varphi=\angle(x, qx)$, then
   \begin{equation*}
   \angle(x, qx)=\angle(qx, q^{2}x)=\angle(x, q^{2}x)=\varphi,\quad \varphi\in \big(0,\dfrac{2\pi}{3}\big);
   \end{equation*}
    \item [(iii)] An orthogonal $q$-basis of $T_{p}M$ exists.
\end{itemize}

\section{The components of the curvature tensor}\label{sec:3}

The Levi-Civita connection on a Riemannian manifold is denoted by $\nabla$. For the Christoffel symbols $\Gamma_{ij}^{s}$ of $\nabla$ it is known that
\begin{equation}\label{gamma}
    2\Gamma_{ik}^{h}=g^{ht}(\partial_{i}g_{tk}+\partial_{k}g_{ti}-\partial_{t}g_{ik}),
    \end{equation}
where $g^{ij}$ are the components of the inverse matrix of $(g_{ij})$.

The curvature tensor $R$ of $\nabla$ is defined by  $$R(x, y)z=\nabla_{x}\nabla_{y}z-\nabla_{y}\nabla_{x}z-\nabla_{[x,y]}z $$  and the local components of $R$ are
\begin{equation}\label{r}
    R_{ijk}^{h}=\partial_{j}\Gamma_{ik}^{h}-\partial_{k}\Gamma_{ij}^{h}+\Gamma_{ik}^{t}\Gamma_{tj}^{h}-\Gamma_{ij}^{t}\Gamma_{tk}^{h}.
\end{equation}
The corresponding tensor $R$ of type $(0, 4)$ is determined as follows
\begin{equation*}
    R(x, y, z, u)=g(R(x, y)z,u).
\end{equation*}

For $(M,g,q)$, we denote $D=(A-B)(A+2B)$ and
\begin{equation}\label{aibi}
  A_{i}=\frac{\partial A}{\partial X^{i}}\ ,\quad B_{i}=\frac{\partial B}{\partial X^{i}}\ ,
\end{equation}
where $A$ and $B$ are the functions from (\ref{f2}).

The inverse matrix of $g$ is
\begin{equation}\label{g-obr}
    (g^{ij})=\frac{1}{D}\begin{pmatrix}
      A+B & -B & -B \\
      -B & A+B & -B \\
      -B & -B & A+B \\
    \end{pmatrix}.
    \end{equation}
Then by direct calculations, having in mind \eqref{f2}, \eqref{gamma}, \eqref{r}, \eqref{aibi} and \eqref{g-obr},
 we obtain
\begin{theorem}\label{thmR}
The nonzero components of the curvature tensor $R$ of type $(0, 4)$ of the manifold $(M, g, q)$ are
\begin{equation*}
\begin{split}
   &R_{1212} = \frac{1}{2}(2B_{21}-A_{11}-A_{22})\\&+\frac{A+B}{4D}\Big(2A_{3}B_{2}-A_{3}^{2}+(B_{1}-B_{2}-B_{3})(B_{1}+B_{2}-B_{3})\Big)\\&-\frac{2B}{4D}\Big((A_{1}-B_{2})(B_{1}+B_{2}-B_{3})-A_{1}A_{3}+A_{3}B_{2}\Big) ,
\end{split}
\end{equation*}
\begin{equation*}
\begin{split}
   &R_{1313} = \frac{1}{2}(2B_{31}-A_{11}-A_{33})\\&+\frac{A+B}{4D}\Big(2A_{2}B_{3}-A_{2}^{2}+(-B_{1}+B_{2}+B_{3})(-B_{1}+B_{2}-B_{3})\Big)\\&-\frac{2B}{4D}\Big((A_{1}-B_{3})(B_{1}-B_{2}+B_{3})-A_{1}A_{2}+A_{2}B_{3}\Big),
\end{split}
\end{equation*}
\begin{equation*}
\begin{split}
   &R_{2323} = \frac{1}{2}(2B_{23}-A_{22}-A_{33})\\&+\frac{A+B}{4D}\Big(2B_{3}A_{1}-A_{1}^{2}+(B_{1}-B_{2}+B_{3})(B_{1}-B_{2}-B_{3})\Big)\\&-\frac{2B}{4D}\Big((A_{2}-B_{3})(B_{2}-B_{1}+B_{3})-A_{1}A_{2}+A_{1}B_{3}\Big),
\end{split}
\end{equation*}
\begin{equation*}
\begin{split}
   &R_{1213} = \frac{1}{2}(B_{21}+B_{31}-B_{11}-A_{23})\\&+\frac{A+B}{4D}\Big(A_{1}(B_{2}-B_{3}+B_{1})+2B_{3}(B_{3}-B_{2}-B_{1})+A_{2}A_{3})\Big)\\&-\frac{B}{4D}\Big(A_{1}^{2}+A_{2}^{2}+A_{3}^{2}+2A_{1}(A_{2}-B_{3})-2A_{2}B_{3}\\&-2A_{3}(B_{1}-B_{3})+(B_{1}-B_{2}-B_{3})(B_{1}+B_{2}-B_{3})\Big),
\end{split}
\end{equation*}
\begin{equation*}
\begin{split}
   &R_{1223} = \frac{1}{2}(B_{22}-B_{12}-B_{23}+A_{13})\\&+\frac{A+B}{4D}\Big(A_{2}(B_{2}+B_{3}-B_{1})-(2B_{3}-A_{1})(2B_{2}-A_{3})\Big)\\&-\frac{B}{4D}\Big(A_{2}^{2}-A_{1}^{2}+A_{3}^{2}+2A_{1}(B_{2}+B_{3})+2A_{2}(B_{2}-B_{3})\\&+2A_{3}(B_{3}-B_{1})-4B_{2}B_{3}+(B_{1}+B_{2}-B_{3})(B_{1}-B_{2}-B_{3})\Big),
\end{split}
\end{equation*}
\begin{equation*}
\begin{split}
   &R_{1323} = \frac{1}{2}(B_{23}-B_{33}+B_{13}-A_{12})\\&+\frac{A+B}{4D}\Big((2B_{2}-A_{1})(2B_{3}-A_{2})-A_{3}(-B_{1}+B_{2}+B_{3}\Big)\\&-\frac{B}{4D}\Big(A_{1}^{2}-A_{2}^{2}-A_{3}^{2}-2A_{1}(B_{2}+B_{3})+2A_{2}(B_{1}-B_{2})\\&+2A_{3}(B_{2}-B_{3})+4B_{2}B_{3}+(-B_{1}+B_{2}+B_{3})(B_{1}-B_{2}+B_{3})\Big).
\end{split}
\end{equation*}
The rest of the nonzero components are obtained by the properties $$R_{ijkh}=R_{khij}, \ R_{ijkh}=-R_{jikh}=-R_{ijhk}\ .$$
\end{theorem}

\section{Some sectional curvatures}\label{sec:4}
In \cite{1}, for $(M, g, q)$ it is proved that $\nabla q=0$ implies
\begin{equation}\label{V1}
  R(x, y, qz, qu)=R(x, y, z, u).
\end{equation}
Therefore it follows the identity
\begin{equation}\label{V2}
  R(qx, qy, qz, qu)=R(x, y, z, u),
\end{equation}
which defines a more general class of manifolds $(M, g, q)$ than the class with the condition $\nabla q=0$.

Let $\{x, y\}$ be a non-degenerate $2$-plane spanned by vectors $x,y\in T_{p}M$, $p\in M$. Then its sectional curvature is
\begin{equation}\label{3.3}
    \mu(x,y)=\frac{R(x, y, x, y)}{g(x, x)g(y, y)-g^{2}(x, y)}.
\end{equation}

\begin{prop}\label{prop1}
Let $(M,g,q)$ be a manifold with property \eqref{V2} and a vector $x$ induce a $q$-basis. Then  $$\mu(x, qx)=\mu(qx, q^{2}x)=\mu(x, q^{2}x).$$
\end{prop}

\begin{proof}
  From (\ref{V2}) we get
\begin{equation}\label{3.6}
    R(q^{2}x, q^{2}y, q^{2}z, q^{2}u)=R(qx, qy, qz, qu)=R(x, y, z, u).
    \end{equation}
    In (\ref{3.6}) we substitute $qx$ for $y$, $x$ for $z$ and $qx$ for $u$
    and we find
    \begin{equation}\label{dopl}
    R(q^{2}x, x, q^{2}x, x)=R(qx, q^{2}x, qx, q^{2}x)=R(x, qx, x, qx).
    \end{equation}
Then, from \eqref{2.12} and \eqref{3.3} it follows \eqref{3.3}.
\end{proof}
Let $x$ induce a $q$-basis of $T_{p}M$ and $\sigma=\{x, qx\}$ be a $2$-plane. It is easy to see that if $y \in\sigma$ and $y\neq x$, then $qy \notin \sigma$. Consequently, $\sigma$ has only two $q$-bases: $\{x, qx\}$ and $\{-x, -qx\}$. That's why the sectional curvature $\mu(x, qx)$ depends only on $\varphi=\angle(x, qx)$.
So, we denote $\mu(x, qx)=\mu(\varphi)$.

\begin{theorem}\label{th4}
Let $(M, g, q)$ be a manifold with property \eqref{V2}. If a vector $u$ induces a $q$-basis, then
\begin{equation}\label{mu-r2}
   \mu(\varphi)=\frac{1-2\cos\varphi}{1+\cos\varphi}\mu\Big(\frac{\pi}{2}\Big)+\frac{3\cos\varphi}{1+\cos\varphi}\mu\Big(\frac{\pi}{3}\Big)\ ,
\end{equation}
where $\varphi=\angle(u, qu)$.
\end{theorem}

\begin{proof}
In (\ref{3.6}) we substitute $qx$ for $y$, $q^{2}x$ for $z$ and $x$ for $u$ and we get
    \begin{equation}\label{dop2}
    R(q^{2}x, x, qx, q^{2}x)=R(qx, q^{2}x, x, qx)=R(x, qx, q^{2}x, x).
    \end{equation}
Let a vector $x$ induce an orthormal $q$-basis.
  If $u=\alpha x+\beta qx +\gamma q^{2}x$, where $\alpha,\beta,\gamma\in \mathbb{R}$, then $qu=\gamma x+\alpha qx +\beta q^{2}x$.
 Due to the linear properties of the curvature tensor $R$, we obtain
 \begin{equation*}
    \begin{split}
    R(u,qu,u,qu)&= (\alpha^{2}-\beta\gamma)^{2}R(x, qx, x, qx)\\&+(\gamma^{2}-\alpha\beta)^{2}R(x, q^{2}x, x, q^{2}x)\\&+(\beta^{2}-\alpha\gamma)^{2}R(qx, q^{2}x, qx, q^{2}x)\\&+2(\alpha^{2}-\beta\gamma)(\gamma^{2}-\alpha\beta)R(x, qx, q^{2}x, x)\\&+2(\gamma^{2}-\alpha\beta)(\beta^{2}-\alpha\gamma)R(q^{2}x, x, qx, q^{2}x)\\&+2(\alpha^{2}-\beta\gamma)(\beta^{2}-\alpha\gamma)R(x, qx, qx, q^{2}x).
    \end{split}
\end{equation*}
Having in mind \eqref{dopl} and \eqref{dop2} we find
 \begin{equation}\label{cor}
    \begin{split}
    &R(u,qu,u,qu)=  \Big((\alpha^{2}-\beta\gamma)^{2}\\&+(\gamma^{2}-\alpha\beta)^{2}+(\beta^{2}-\alpha\gamma)^{2}\Big)R(x, qx, x, qx)\\&+2\Big((\alpha^{2}-\beta\gamma)(\gamma^{2}-\alpha\beta)+(\gamma^{2}-\alpha\beta)(\beta^{2}-\alpha\gamma)\\&+(\alpha^{2}-\beta\gamma)(\beta^{2}-\alpha\gamma)\Big)R(x, qx, q^{2}x, x).
    \end{split}
\end{equation}

Since $\{x, qx, q^{2}x\}$ is an orthonormal $q$-basis, we have
\begin{equation*}
    g(u, u)= g(qu, qu)=\alpha^{2}+\beta^{2}+\gamma^{2},\quad g(u, qu)= \alpha\beta+\beta\gamma+\gamma\alpha.
\end{equation*}
     We suppose that $g(u, u)=1$. From \eqref{2.12} and \eqref{3.3} we get
\begin{equation}\label{mu}
    \mu(\varphi)= \frac{R(u,qu,u,qu)}{1-
    \cos^{2}\varphi}\ ,
\end{equation}
and \begin{equation*}
   \alpha^{2}+\beta^{2}+\gamma^{2}=1,\quad
   \alpha\beta+\beta\gamma+\gamma\alpha=\cos\varphi.
\end{equation*}
We express $\alpha, \beta, \gamma$ by $\cos\varphi$ as follows:
    \begin{align*}
    (\cos\varphi)^{2}-\cos\varphi &=(\alpha^{2}-\beta\gamma)(\gamma^{2}-\alpha\beta)+(\gamma^{2}-\alpha\beta)(\beta^{2}-\alpha\gamma)\\&+(\alpha^{2}-\beta\gamma)(\beta^{2}-\alpha\gamma), \\
      1-(\cos\varphi)^{2}&=(\alpha^{2}-\beta\gamma)^{2}+(\gamma^{2}-\alpha\beta)^{2}+(\beta^{2}-\alpha\gamma)^{2}.
    \end{align*}
Then, from \eqref{cor} and \eqref{mu}, we obtain
\begin{equation}\label{mu-r}
   \mu(\varphi)=\mu\Big(\frac{\pi}{2}\Big)+\frac{2\cos \varphi}{1+\cos\varphi}R(x,qx,x,q^{2}x).
\end{equation}
In \eqref{mu-r} we substitute $\dfrac{\pi}{3}$ for $\varphi$ and we get
    $$R(x, qx, x, q^{2}x)=\dfrac{3}{2}\Big(\mu\Big(\dfrac{\pi}{3}\Big)-\mu\Big(\dfrac{\pi}{2}\Big)\Big).$$
The last result and \eqref{mu-r} imply \eqref{mu-r2}.
\end{proof}
\begin{corollary}\label{th5}
Let $(M, g, q)$ be a manifold with property \eqref{V1}.
If a vector $u$ induces a $q$-basis, then
\begin{equation}\label{mu-r3}
   \mu(\varphi)=\frac{1-\cos\varphi}{1+\cos\varphi}\mu\Big(\frac{\pi}{2}\Big)\ ,
\end{equation}
where $\varphi=\angle(u, qu)$.
\end{corollary}

\begin{proof}
 Since \eqref{V1} is valid, we find
 \begin{equation}\label{V13}
 R(x, y, z, u)=R(x, y, qz, qu)= R(x, y, q^{2}z, q^{2}u).
\end{equation}
In \eqref{V13} we substitute $qx$ for $y$, $x$ for $z$ and $qx$ for $u$ and we get
         \begin{equation}\label{dopv1}
    R(x, qx, x, q^{2}x)=-R(x, qx, x, qx).
    \end{equation}
If we suppose that $\{x, qx, q^{2}x\}$ is an orthonormal $q$-basis, then from \eqref{mu-r} and \eqref{dopv1} follows \eqref{mu-r3}.
\end{proof}

\section{An example of $(M, g, q)$}\label{sec:5}

In \cite{1} it is proved that, the structure $q$ is parallel with respect to the Levi-Civita connection $\nabla$ of $g$ on a manifold $(M, g, q)$ if and only if the gradients of $A$ and $B$ satisfy the following equality:
\begin{equation}\label{parallel}
    \grad A=\grad B\begin{pmatrix}
      -1 & 1 & 1 \\
      1 & -1 & 1 \\
      1 & 1 & -1 \\
    \end{pmatrix}.
\end{equation}

In this section we discuss an example of a manifold $(M, g, q)$ which satisfies \eqref{V2}, but doesn't  satisfy \eqref{parallel}.

\begin{theorem}\label{prop}
The property \eqref{V2} of the manifold $(M, g, q)$ is equivalent to the conditions
\begin{equation}\label{r1=r6}
    R_{1212}=R_{1313}=R_{2323},\quad R_{1213}=R_{1323}=-R_{1223},
\end{equation}
where $R_{ijkh}$ are the local components of the curvature tensor $R$ of type $(0, 4)$.
\end{theorem}

\begin{proof}
The local form of \eqref{V2} is
\begin{equation}\label{rlocal}
    R_{tslm}q_{i}^{t}q_{j}^{s}q_{k}^{l}q_{h}^{m}=R_{ijkh}.
\end{equation}
From \eqref{f2} and \eqref{rlocal} we find
\begin{align}\label{r1=r3}\nonumber
    R_{1212}=R_{2323},\quad R_{1313}=R_{2121},\\ R_{2321}=R_{1213},\quad
    R_{2331}=R_{1223},\\\nonumber R_{2131}=R_{1323},\quad R_{3131}=R_{2323},
\end{align}
which implies \eqref{r1=r6}.

Vice versa, from \eqref{r1=r6} it follows \eqref{r1=r3}. Having in mind \eqref{f2} we get \eqref{rlocal}.
\end{proof}

Let $(M, g, q)$ be a manifold with
\begin{equation}\label{36}
    A=2X^{1}\ ,\ B=2X^{1}+X^{2}+X^{3} ,
\end{equation}
where $$2X^{1}+X^{2}+X^{3}>0,\quad X^{2}+X^{3}<0.$$
Obviously, the condition \eqref{f3} is satisfied.
Due to \eqref{f2}, \eqref{g-obr}, \eqref{36} and  Theorem~\ref{thmR} we obtain
   \begin{equation}\label{comR}
   \begin{split}
       R_{1212} &=R_{1313}=R_{2323}= -\frac{B}{(A-B)(A+2B)}\ ,\\
   R_{1213} &=R_{1323}=R_{1223}=0 .
    \end{split}
   \end{equation}
We check directly that the conditions \eqref{r1=r6} are valid, but the conditions \eqref{parallel} for the functions $A$ and $B$ are not valid.

Consequently, we obtain the following
\begin{theorem}
The manifold $(M,g, q)$ with \eqref{36} satisfies the curvature identity \eqref{V2}.
Furthermore, the structure $q$ is not parallel with respect to the Levi-Civita connection $\nabla$ of $g$.
\end{theorem}
The Ricci tensor $\rho$ and the scalar curvature $\tau$ are given by the well-known formulas:
\begin{equation}\label{def-rho}
    \rho(y,z)=g^{ij}R(e_{i}, y, z, e_{j}),\
    \tau=g^{ij}\rho(e_{i}, e_{j}).
\end{equation}
We obtain the components of $\rho$ and the value of $\tau$:
\begin{equation}\label{rho32}
\begin{split}
    \rho_{12}&=\rho_{13}=\rho_{23}=\frac{B^{2}}{(A-B)^{2}(A+2B)^{2}}\ ,\\\rho_{11}&=\rho_{22}=\rho_{33}=\frac{2B(A+B)}{(A-B)^{2}(A+2B)^{2}}\ ,
  \end{split}
\end{equation}
\begin{equation}\label{tau32}
    \tau=\frac{6AB}{(A-B)^{3}(A+2B)^{2}}\ .
\end{equation}

Therefore, we arrive at the following
\begin{prop}\label{kt4}
For the manifold $(M, g, q)$ with \eqref{36}, the following assertions are valid:
\begin{itemize}
\item[(i)]  The components of the curvature tensor $R$  are \eqref{comR}, i.e. $M$ is not a flat manifold;

\item[(ii)]  The components of the Ricci tensor $\rho$ are \eqref{rho32};

\item[(iii)] The scalar curvature $\tau$ is \eqref{tau32}.
\end{itemize}
\end{prop}



\end{document}